\documentclass[a4paper,oneside,reqno,11pt]{amsart}
\usepackage[english]{babel}
\usepackage{amsfonts}
\usepackage{amsmath}
\usepackage{amssymb}
\usepackage{amsthm}
\usepackage{enumerate}
\usepackage{geometry}

\title
[Convergence of solutions of SDEs to Harris flows]
{Convergence of solutions of SDEs to Harris flows}
\author{M.~B.~Vovchanskii}
\address{M.~B.~Vovchanskii: Department of Theory of Random Processes, Institute of Mathematics, National Academy of Sciences of Ukraine, Tereshchenkivska Str.~3, Kiev~01601, Ukraine}
\email{vovchansky.m@gmail.com}

\usepackage{amsmath, amsthm, amsfonts, amssymb, bbm}
\usepackage[english]{babel}
\usepackage{mathtools}

\renewcommand{\lg}{\langle}
\newcommand{\rg}{\rangle}

\newcommand{\vf}{\varphi}

\newcommand{\1}{\mathbbm{1}}
\newcommand{\mbR}{{\mathbb R}}
\newcommand{\mbN}{{\mathbb N}}

\newcommand{\mcD}{\mathcal{D}}
\newcommand{\mcM}{\mathcal{M}}
\newcommand{\mcC}{\mathcal{C}}
\newcommand{\mcI}{\mathcal{I}}
\newcommand{\mcP}{\mathcal{P}}

\newcommand{\mfA}{\mathcal{A}}

\newcommand{\e}{\mathrm{e}}
\newcommand{\Cov}{\mathrm{Cov}}
\newcommand{\E}{\mathrm{E}}
\newcommand{\ve}{\varepsilon}

\DeclareMathOperator\supp{supp}

\newtheorem{theorem}{Theorem}

\newtheorem{proposition}{Proposition}
\newtheorem{definition}{Definition}
\newtheorem{remark}{Remark}

\begin{document}

\begin{abstract}
A method of the approximation of a coalescing Harris flow with homeomorphic stochastic flows built as solutions to SDEs w.r.t. continuous martingales with spatial parameters in the sense of Kunita is proposed. The joint convergence of forward and backward flows as diffusions is obtained, as well as the joint convergence of forward and backward transformations of the real axe under the action of the flows.  
\end{abstract}

\maketitle

We consider a system of correlated and sticking together after a meeting Brownian motions with $\mbR\times\mbR^{+}$ as a set of start points. The object discussed is characterized via its $n-$point motions that are approximated with $n-$point motions of ''ordinary'' stochastic flows, where hereinafter a term ''stochastic flow'' is referring to a family of random mappings. The first class of examples considered in the paper is provided with flows of solutions to SDEs w.r.t. a martingale with spatial parameters in the sense of~\cite{kunita}, while the second one is provided with Harris flows~\cite{harris}. 

The Harris flow's original definition is up to Harris himself~\cite{harris}, though we use a modified formulation (see~\cite{dorogovtsev:mono, warren.watanabe}). Let $\mcD(\mbR)$ be a separable topological space of rcll functions on $\mbR$ equipped with the Skorokhod topology~\cite{vogel, billingsley}. The space $\mcD(\mbR)$ is completely metrizable (see \cite{billingsley} for an example of such a metric). 
A composition $f(g)$ is denoted $f\circ g,$ and $Id$ is the identity mapping. 
\begin{definition}
\label{harris.flow}
A Harris flow $X$ with the infinitesimal covariance $\vf$ is a family of $\mcD(\mbR)-$valued random variables $\{X(s,t)\mid X(s,t)\equiv X(\cdot, s, t), s\le t\}$ such that
\begin{enumerate}
\item
for any $s\le t\le r$ $P\left\{X(\cdot,s,r)=X(\cdot,t,r)\circ X(\cdot,s,t)\right\}=1$; $X(s,s)=Id$ a.s.;
\item
for any $t_1\le t_2 \le \ldots \le t_n$ random elements $X(t_1,t_2),\ldots,X(t_{n-1},t_n)$ are independent;
\item
for any $s, t\in\mbR, h>0$ $Law\left(X(s,t)\right)=Law\left(X(s+h,t+h)\right);$
\item
as $h\rightarrow 0+,$ $X(0,h)\rightarrow Id$ in probability;
\item
for any $x$ a process $t\mapsto X(x,0,t)-x$ is a Brownian motion started at $0$ w.r.t. filtration $\sigma\left\{ X(u_1,u_2), 0\le u_1\le u_2 \le t\right\}_{t\ge 0};$
\item
for any $x,y$ $\left\langle X(x,0,\cdot), X(y, 0,\cdot)\right\rangle(t) = \int_0^t \vf\left(X(x, 0,s)-X(y,0,s)\right)ds.$
\end{enumerate} 
\end{definition}
	
It is shown in~\cite{harris} that if a symmetric continuous nonnegative definite function $\varphi$ whose Fourier transform is not of pure jump type is Lipschitz continuous outside any interval containing the origin the corresponding Harris flow exists. If $\vf$ is smooth enough, the flow can be considered as a flow in the sense of Kunita~\cite{kunita, warren.watanabe}. However, not every Harris flow is a flow of homeomorphisms in contrast to those treated in~\cite{kunita}. Indeed, the difference $\xi=\left(X(u_1,0,\cdot)-X(u_2,0,\cdot)\right)$ is a Feller diffusion with the infinitesimal operator $\left(1-\varphi(x)\right)\frac{d^2}{dx^2}$ on the upper halfline until it hits the origin (if ever). If $\int_0^\delta \frac{x}{1-\varphi(x)}dx$ is finite for some small $\delta$ the Feller criterion~\cite{cherny.engl, harris} implies that $\xi$ hits the origin in a finite time a.s.. It being a nonnegative martingale, the process $\xi$ never leaves $0$ after hitting it. Thus Harris flows can be referred to as coalescing flows in this case. At the same time, the origin can still be either an exit or a regular point for the diffusion $\xi$ (see \cite{tsirelson.products}[Section 9c]). In \cite{warren.watanabe} it is shown how the regularity of the origin influences properties of the noises associated with coalescing Harris flows.

The Brownian web (\cite{dynamics, BW:full, BW:convergence, selfrepelling}) can be considered as an extreme example of a Harris flow with a discontinuous infinitesimal characteristic $\vf(x)=\1_{x=0}.$ In~\cite{dorogovtsev:brownian.flow} a constructive approach to $n-$point motions of the Brownian web is proposed that is similar to the one adopted in this paper. More specifically,
let $W$ be a Brownian sheet~\cite{kotelenez}. Suppose that a sequence $\{\kappa_n\}_{n\ge 1}$ of infinitely differentiable functions weakly converges to $\delta_0$ in the sence of generalized functions, and there exist infinitely differentiable square integrable functions $\rho_n$ such that $\rho_n\ast\rho_n=\kappa_n, n\in\mbN.$ One consider the following SDE
\begin{equation}
\label{int.brownian.sheet}
X_n(x,t) = x + \int_0^t\int_{\mbR} \rho_n(y-X_n(x,s))W(ds,dy),
\end{equation}
(for a general theory of integration w.r.t. the Brownian sheet, the reader is referred to~\cite{kotelenez, dorogovtsev:mono}). Given $x_1,\ldots,x_n\in\mbR$ a sequence $\left\{ X_n(x_1,\cdot),\ldots, X_n(x_n,\cdot) \right\}_{n\ge 1}$  weakly converges to $\left(X(x_1,0,\cdot),\ldots, X(x_n,0,\cdot) \right)$ in $\left(\mcC(\mbR^+)\right)^n,$  $X$ being a Brownian web~\cite{dorogovtsev:brownian.flow}. 

In this paper we obtain a similar result for a Harris flow whose infinitesimal covariance is a characteristic function of a symmetric stable law, possibly, with infinite mean. In the latter case the convolutional square root of the covariance does not belong to $L^2(\mbR),$  violating conditions for~\eqref{int.brownian.sheet} to have a solution. 
In order to overcome this difficulty we use integration w.r.t. continuous martingales with spatial parameters~\cite{kunita} obtaining $n-$point motions of a coalescing Harris flow as weak limits of $n-$point motions of flows of homeomorphisms given by the SDE
\begin{equation}
\label{int:kunita.mart}
X_\ve(x, s, t) = x + \int_s^t F_{\ve}(X_\ve(x,s,r), dr),
\end{equation}
where $F_\ve$ is a continuous $\mcC(\mbR)-$valued martingale as defined in~\cite{kunita}. Now we give a rigorous description of the approach proposed.

We consider 
\begin{equation*}
\vf(x) = \e^{-\beta|x|^{\alpha}}, x \in \mbR, \beta\in(0;+\infty),\alpha\in(0;2).
\end{equation*} 
Suppose that $\{\vf_\ve\}_{\ve\in(0;1)}$ is a fixed throughout the paper sequence of twice continuously differentiable symmetric nonnegative definite functions such that $\vf_\ve\to \vf, \ve\to0+,$ uniformly on compact subsets of $\mbR,$ and $\vf_\ve(0)=1.$ For instance, one can consider a mollifier $h$ and put
\begin{equation}
\label{smoothing}
\vf_{\ve} =  \ve^{-1} c_\ve \vf \ast h\left(\frac{\cdot}{\ve}\right), \ve \in (0;1),
\end{equation}
with $c_\ve$ selected to give $\vf_\ve(0)=1.$ 
A reference example of such mollifier is provided with a Gaussian density. In this case the result of~\cite{harris}[Lemma 10.4] states that $n-$point motions of Harris flows with the infinitesimal covariances built as in~(\ref{smoothing}) weakly converge to those of a Harris flow with the infinitesimal covariance $\vf.$ In the paper we investigate the joint convergence of $n-$point motions of forward and backward (inverse) flows viewed as diffusions (see \cite{kunita}[Chapter 5]) and the convergence as $\mcD(\mbR)-$valued random elements of Definition~\ref{harris.flow}. Note that in the case of a smooth infinitesimal covariance the inverse flow $X^{-1}$ solves an SDE analogous to that of~\eqref{int:kunita.mart}, although in the inverse time and w.r.t a backward infinitesimal generator~\cite{kunita} so $X^{-1}(\cdot, s, t)$ is a homeomorphism of $\mbR$ onto itself. Moreover, since the infinitesimal covariance is symmetric additionally the flows $X^{-1}$ and $X$ have the same distribution~\cite{kunita}[Theorem 4.2.10]. In constrast, for $\vf$ in question a mapping $X(\cdot, s,t)$ is a.s. a step function~\cite{harris, matsumoto}. However, one can still consider an inverse flow~\cite{harris, korenovskaya.ukr.mat} on $[t_1;t_2]$ defined via 
\begin{align}
\label{inverse.flow:definition}
X^{-1}(x,t_1,t_2,s)  &= \inf\big\{X(y,r,t_1+t_2-s)\mid X(y,r,t_2)\ge x, y\in\mbR, r \in[t_1;t_1+t_2-s]\big\}, \nonumber \\
& s \in[t_1;t_2],
\end{align}
that is, the inverse flow is built with using all possible trajectories of the forward flow on $[t_1;t_2].$ In the case of homeomorphic flows this definition of the inverse flow coincides with the standard one formulated in terms of integrals w.r.t. a backward infinitesimal generator~\cite{kunita}, due to the ordered nature of the trajectories within the flow. But for flows of discontinuous mappings we understand this relations as a definition of the inverse flow (see also~\cite{BW:full,BW:convergence, selfrepelling}). It is worth noting that for any Harris flow $X$ mappings $X^{-1}(\cdot,0,t,s)$ and $X(\cdot,s,t)$ have the same distribution~\cite{harris}[Theorem 10.5] as elements of $\mcD(\mbR).$ 
		
The following notation is adopted hereinafter. Let $\mcC(\Delta),$ where $\Delta$ can be a finite closed interval, the real line or the real halfline, be a space of continuous functions on $\Delta$ equipped with the topology of uniform convergence on compact sets. In product spaces the product topology is always considered henceforth. Denote $\mcC^\infty([s;t]) = \left(\mcC([s;t])\right)^\infty, \mcC^N([s;t]) = \left(\mcC([s;t])\right)^N, N\in\mbN.$  Given real numbers $a, a_1, b\colon a < a_1 < b$ a function $f\in\mcC([a_1;b])$  can be embedded into $\mcC([a;b])$ with a transformation $\mcP_{a,b}f(s)=\1_{s\in[a;a_1]}f(a_1)+f(s)\1_{s\in(a_1;b]}, s\in[a;b].$
\begin{remark}
Given nonnegative $s,r,t\colon s <r<t$ and an arbitrary Harris flow $X$ one has $\mcP_{s,t}X(x,r,\cdot)\in\mcC([s;t])$ and $\mcP_{s,t}X^{-1}(x,r,t,r+t-\cdot)\in\mcC([s;t]).$ Note that $\mcP_{s,t}X^{-1}(x,r,t,r+t-\cdot)(u) = X^{-1}(x,r,t,t), u \in[s;r].$ 
\end{remark}

We start with collections of $\mcC(\mbR)-$valued Brownian motions $F_\ve, \ve\in(0;1)$~\cite{kunita}.
	
\begin{proposition}
 \label{gaussian.process.existence}
Given $\ve\in(0;1)$ there exists a real-valued Gaussian process $F_\ve\equiv\{F_{\ve}(x, t)\mid x\in\mbR, t\in\mbR^{+}\}$ such that 
 \begin{enumerate}
 \item $\forall t\in\mbR^+ \quad F_{\ve}(\cdot, t) \in \mcC(\mbR);$
  \item $\forall t_1 < \ldots < t_n$ $F_{\ve}(\cdot, t_1), F_{\ve}(\cdot, t_2) - F_{\ve}(\cdot, t_1), \ldots, F_{\ve}(\cdot, t_n) - F_{\ve}(\cdot, t_{n-1})$ are independent $\mcC(\mbR)-$valued elements;
  \item $\forall t\ge 0, s\ge0 \ F_{\ve}(\cdot, t + s) - F_{\ve}(\cdot, s) \overset{d}{=} F_{\ve}(\cdot, t) \overset{d}{=} \sqrt{t} F_{\ve}(\cdot, 1);$
  \item the process $t\to F_{\ve}(\cdot, t)\in \mcC(\mbR)$ is a.s. continuous;
  \item $\Cov(F_{\ve}(t,x), F_{\ve}(s,y)) = \min\{t, s\}\vf_{\ve}(x-y).$
 \end{enumerate}
\end{proposition}
\proof Since the mapping $(t,s,x,y)\to \min\{t, s\} \cdot \vf_{\ve}(x-y)$ is nonnegative definite as a product of two covariances a Gaussian process $F_\ve$ with such covariance exists. To check the continuity of $F_\ve$ one calculates     
\begin{align*}
 \rho_T(u) & \coloneqq \sup_{(x-y)^2 + (t-s)^2 \le u^2, t,s\le T}\sqrt{\E(F_\ve(x,t) - F_\ve(y, s))^2}  \\
 & = \sup_{(x-y)^2 + (t-s)^2 \le u^2, t,s\le T}\sqrt{ t + s - 2\left( t\wedge s \right) \vf_\ve(x-y)} \\
 & =\sup_{(x-y)^2 + (t-s)^2 \le u^2, t,s\le T}\sqrt{(t + s) (1 - \vf_{\ve}(x-y)) + (t\vee s-t\wedge s)\varphi_\ve(x-y) }  \\
& \le \sup_{(x-y)^2 + (t-s)^2 \le u^2, t,s\le T}\sqrt{2T |\vf_\ve''(0)| u^2 + u} \le C \sqrt{T} u,
\end{align*}
with some constant $C$. 
Since for sufficiently small $\delta$
$$
\int_\delta^\infty \rho_T(e^{-x^2}) dx<\infty 
$$
the process $\{F_\ve(x,t)\mid x\in\mbR, t\in[0;T]\}$ has a continuous version w.r.t. both arguments~\cite{adler.taylor}. Applying an usual expansion argument we get the existence of such version for $\{F(x,t)\mid x\in\mbR, t\in\mbR^+\}.$ 
This is the version we work with. Property 4 is checked trivially. 

Since $F_\ve$ is a Gaussian field  Properties 2 and 3 immediately follow from usual calculations of corresponding covariances.
\qed	

The process $F_\ve$ is a continuous $\mcC(\mbR)-$valued martingale in the sense of Kunita.

\begin{proposition}
\label{harris.flow.prelimit.existence}
Fix $\ve\in(0;1).$ Then there exists a Harris flow $\{X_\ve(\cdot, s, t)\mid 0\le s \le t\}$ with the infinitesimal covariance $\vf_\ve$ such that 
\begin{enumerate}
\item
for any $x\in\mbR, 0\le s \le t$
\begin{equation*}
\label{main.sde}
X_\ve(x, s, t) = x + \int_s^t F_{\ve}(X_\ve(x,s,r), dr);
\end{equation*}
\item
$\forall 0\le s\le t$ $X_\ve(\cdot, s, t)$ is a homeomorphism on $\mbR.$
\end{enumerate}
\end{proposition}
Proof. We use Theorem 4.5.1 from~\cite{kunita} which states the existence of such a flow on a finite time interval. To be applicable, Theorem 4.5.1 imposes additional requirements on $\vf_\ve$ that in our case are reduced to the finiteness of
$$
\sup_{x,x^\prime, y, y^\prime\in K, x\not=x^\prime, y\not=x^\prime} \frac{|\vf_\ve(x-y)+\vf_\ve(x^\prime-y^\prime)-\vf_\ve(x-y^\prime)-\vf_\ve(x^\prime-y) |}{|x-x^\prime||y-y^\prime|}
$$
for any compact subset $K$ of $\mbR,$ which is a consequence of the mean value theorem since $\vf_\ve$ has bounded second partial derivatives. 
\qed
	
Let $\mcM(\mbR)$ be a space of locally finite nonnegative Radon measures on the real line equipped with the vague topology. Put $\mcM^N(\mbR) = \left(\mcM(\mbR)\right)^N,$ $N\in\mbN.$ The space $\mcM^N(\mbR)$ is separable~\cite{kallenberg:random.measures}.
\begin{theorem}
\label{harris:approximation}
Let $\{X_\ve\}_{\ve\in(0;1)}$ be the Harris flows from Proposition~\ref{harris.flow.prelimit.existence}, and let $X$ be a Harris flow with the infinitesimal covariance $\vf.$ Fix $T>0$ and a set $\{(x_n,t_n)\}_{n\in\mbN}\in\left(\mbR\times[0;T]\right)^{\infty}.$ Then 
\begin{align*}
& \Big(  \mcP_{0,T}X_\ve(x_1, t_1, \cdot), \mcP_{0,T}X^{-1}_\ve(x_1,  t_1, T, T+t_1-\cdot), \ldots,  \\
& \phantom{abcabc} \mcP_{0,T}X_\ve(x_N, t_N, \cdot), \mcP_{0,T}X^{-1}_\ve(x_N, t_N, T, T+t_N-\cdot), \ldots\Big)  \\
& \Rightarrow \Big(  \mcP_{0,T}X(x_1, t_1, \cdot), \mcP_{0,T}X^{-1}(x_1,  t_1, T, T+t_1-\cdot),  \ldots, \\ 
& \phantom{abcabc}  \mcP_{0,T}X(x_N, t_N, \cdot), \mcP_{0,T}X^{-1}(x_N,  t_N, T, T+t_N-\cdot), \ldots \Big), 
\end{align*}
in $\mcC^{\infty}([0;T])$ as $\ve\rightarrow 0+.$

Let $\lambda$ be the Lebesque measure on the real line. For $\ve\in(0;1), 0\le s\le t \le T,$ define the following  $\mcM(\mbR)-$valued random elements:
\begin{align*}
\mu_\ve(s,t) & =  \lambda\circ X_\ve(\cdot, s,t)^{-1}, \\
\mu(s,t) & =  \lambda\circ X(\cdot, s,t)^{-1}, \\
\widehat\mu_\ve(s,t) & =  \lambda\circ \left(X_\ve^{-1}(\cdot,0,t,s)\right)^{-1},\\
 \widehat\mu(s,t) & =  \lambda\circ \left(X^{-1}(\cdot,0,t,s)\right)^{-1}. 
\end{align*}
Then
for any $s_1 \le, \ldots \le s_N, t_1 \le, \ldots \le t_N, s_i \le  t_i, i=\overline{1,N}, N\in\mbN,$
\begin{align*}
& \left(  \mu_\ve(s_1, t_1),\ldots, \mu_\ve(s_N, t_N), \widehat\mu_\ve(s_1, t_1),\ldots, \widehat\mu_\ve(s_N, t_N) \right) \\
& \Rightarrow \left(  \mu(s_1, t_1),\ldots, \mu(s_N, t_N), \widehat\mu(s_1, t_1),\ldots, \widehat\mu(s_N, t_N) \right),
\end{align*}
in $\mcM^{2N}(\mbR)$ as $\ve\rightarrow 0+.$
\end{theorem}
Proof. Ideas and techniques from~\cite{harris, dorogovtsev:large.deviations, portenko.shurenkov} are used in the proof. 

The set $\{(x_n,t_n)\}_{n\in\mbN}$ is additionally supposed to contain all duadic numbers in $\mbR\times[0;T],$ which is always achievable.   

We start with a result on the finite-dimensional convergence. As it was noted previously, \cite{harris}[Lemma 10.4] establishes the convergence of 
$$ 
\left(  X_\ve(x_1, s, \cdot),\ldots, X_\ve(x_N, s, \cdot)\right)
$$ for a particular choice of $\{\vf_\ve\}_{\ve\in(0;1)},$ although a part of the reasoning is left to a reader. Due to our case being more general and in order to keep the presentation comprehensive we present a complete proof covering and generalizing that of~\cite{harris}[Lemma 10.4]. 

Fix a natural number $K.$ Suppose that functions $a_{ij}, b_i, i,j=\overline{1,K},$ are continuous and bounded, a matrix $\|a_{ij}\|_{i,j=\overline{1,K}}$ is nonnegative definite, and define an operator $\mfA$ acting on the space of infinitely differentiable functions with bounded derivatives via
$$
\mfA = \frac{1}{2}\sum_{i,j=\overline{1,K}} a_{ij}\frac{\partial^2}{\partial x_i \partial x_j} +\sum_{i=\overline{1,K}} b_i \frac{\partial}{\partial x_i}.
$$
Let $C$ be a set of continuous functions whose coordinates stay equal after the moment they meet. Fix $s\ge0.$ Given $y\in\mbR^K$ and $s\in\mbR^+$ a measure $P_{y,s}^K$ on $\mcC^K([s;+\infty))$ is called a $C-$solution for the martingale problem for the operator $\mfA$ if for any compactly supported infinitely differentiable $f$ a process $[s,+\infty)\ni r\mapsto f(\omega(r))-\int_s^r \mfA f(\omega(q))dq$ is a martingale w.r.t. $P_{y,s}^K$, where $\omega$ stands for a canonical realisation of a diffusion process, and additionally $P_{y,s}^K\{\omega(s)=y\}=1,$ $\left\{\omega(r+\cdot)\mid \omega \in C \right\}\subset C, $ and $P_{y,s}^K(C)=1.$ 

Solutions $P_{y, s}^K, y\in\mbR^K, s\in\mbR+,$ defined and investigated in~\cite{harris}[Sections 2-3], posess a strong Markov property and are Feller processes. They exist for $\mfA=\frac{1}{2}\sum_{i,j=1}^K g(x_i-x_j)\frac{\partial^2}{\partial x_i \partial x_j},$ $g =\vf_\ve$ or $g=\vf$ and are unique (idem). So since the coefficients of the operator $\mfA$ do not depend on the time variable, all measures $P_{y, s}^K, s \ge 0,$ are essentially shifts of $P_{y,0}^K.$ 
Speaking informally, $\{P^K_{y, 0}\}$ define a process that solves the martingale problem for $\mfA$ up to the moment of hitting the boundary of $\{z\in\mbR^K\mid z_1, \ldots, z_K \mathrm{\ are\ distinct }\},$ after that stays on the boundary and its distinct coordinates again form a solution to the martingale problem for an operator of the same form in the space of a lesser dimension until a new collision happens and so on. In our case each coordinate itself is a Brownian motion. 

Let $s\ge 0, y_1 \le \ldots \le y_K$ be fixed. For any $\ve$ the distribution of 
$$
\eta^K_\ve(y_1,\ldots,y_K)=\left(X_\ve(y_1, s, \cdot),\ldots, X_\ve(y_K, s, \cdot)\right)
$$
is a $C-$solution for the martingale problem for the operator $\frac{1}{2}\sum_{i,j=1}^K \vf_\ve(x_i-x_j)\frac{\partial^2}{\partial x_i \partial x_j}.$ The sequence $\{\eta^K_\ve(y_1,\ldots,y_K)\}_\ve$ is tight in $\mcC^K([s;+\infty))$ by the Tychonoff theorem. Suppose that $\eta^K$ is a weak limit of this sequence as $\ve\rightarrow 0+.$ Note that $\eta^K\in C.$ For this we argue differently to~\cite{harris}[Lemma 10.4]. Each coordinate of $\eta^K$ is a martingale w.r.t. the joint filtration, inheriting this property from the prelimit processes. 
Since $\{(f_{1},\ldots,f_N)\mid f_i-f_j\mathrm{\ never\ changes\ the\ sign\ }\}$ is a closed set in $\mcC^K([s;\infty))$ any difference $\eta^K_k-\eta^K_i$ is a martingale that does not change the sign, which implies the claim.
 
The uniform convergence on compact sets of the sequence $\{\vf_\ve\}_\ve$ to $\vf$ is used to check the represention of the join characteristic of the coordinates of $\eta^K$ on any finite time interval and to prove that the process $r\mapsto f(\eta^K(r))-\int_s^r \mfA f(\eta^K(q))dq$ stays a martingale for $\mfA=\frac{1}{2}\sum_{i,j=1}^K \vf(x_i-x_j)\frac{\partial^2}{\partial x_i \partial x_j}$ (see \cite{portenko.shurenkov}[Remark after Theorem 9 in Chapter 3] for a standard proof). 
Using the same reasoning as that of~\cite{harris}[Lemma 3.2] one gets $Law(\eta^K)=P^K_{y,s}.$

Now we consider the case of distinct moments of start, $s_1 \le \ldots \le s_K.$ Fix $N\in\mbN.$ The sequence 
$$
\xi_\ve=\left(\mcP_{0,T}X_\ve(x_1, s_1, \cdot),\ldots, \mcP_{0,T}X_\ve(x_N, s_N, \cdot) \right)
$$ 
is tight in  $\mcC^N([0;+\infty))$ by the Tychonoff theorem as its each coordinate is a Brownian motion after a certain moment of time and a constant before that. We shall show that the distribution of a weak limit of $\{\xi_\ve\}_{\ve}$  is uniquely determined. Let $\xi=(\xi_1,\ldots,\xi_N)$ be such a limit.

To simplify the notation we suppose that all $s_1,\ldots, s_N$ are distinct; adjustments needed in a general case can be  made easily.  
Fix $M\ge N$ and $r_1,\ldots,r_M$ such that $s_1 \le  r_1 \le \ldots \le r_M.$ Let $\{r_1,\ldots,r_M\}=\bigcup_{i=\overline{1,N}}A_i, A_i \subset [s_i;s_{i+1}), i =\overline{1,N-1}, A_N\subset[s_N;+\infty).$ Note that $\xi_j(r)=\xi_j(s_j), r\le s_j, j= \overline{1,N}.$ Let $f_j,j=\overline{1,M},$ be bounded continuous functions from $\mbR^N$ to $\mbR.$ 
Define the following measurable functions 
\begin{align*}
&g_N(Y) = \int_\Omega \prod_{j\in A_N}f_j(\omega(r_j)) P^N_{Y,s_N}(d\omega), \ Y\in\mbR^N, \\
&g_l(Y) = \int_\Omega \prod_{j\in A_{l}}f_j\big(\omega(r_j)\big) g_{l+1}\big(\omega(s_{l+1})\big) Q^l_{Y}(d\omega), \ Y\in\mbR^N, \\
&Q^l_{Y}=P^l_{(Y_1,\ldots,Y_l),s_l}\otimes \bigotimes_{j=\overline{l+1,N}} \delta(Y_j), \ l=\overline{1,N-1}, 
\end{align*}
where a measure $\delta(u)$ is an atomic measure concentrated at a function identically equal to $u.$ Here the product over an empty set equals $1$ by definition.
Then
\begin{align*}
E&\prod_{j=\overline{1,M}}f_j(\xi(r_j)) = E E\left( \prod_{j=\overline{1,M}}f_j(\xi(r_j)) \mid \xi(r), r \le s_N\right) \\ 
& = E \prod_{j\in\bigcup_{j=\overline{1,N-1}}A_j}f_j(\xi(r_j)) E\left( \prod_{k\in A_N}f_k(\xi(r_k)-\xi(s_N)+\xi(s_N)) \mid \xi(r), r \le s_N\right) \\ 
& = E \prod_{j\in\bigcup_{j=\overline{1,N-1}}A_j}f_j(\xi(r_j)) \left( \int_\Omega \prod_{k\in A_N}f_k\left(\omega(r_k)\right) P^N_{Y,s_N}(d\omega) \right)\Big|_{Y = \xi(s_N)} \\
& = E \prod_{j\in\bigcup_{j=\overline{1,N-1}}A_j}f_j(\xi(r_j)) g_N(\xi(s_N))\\
& = E \prod_{j\in\bigcup_{j=\overline{1,N-2}}A_j}f_j(\xi(r_j)) E\Big(\prod_{i\in A_{N-1}}f_j(\xi(r_i)) g_N(\xi(s_N)) \mid\xi_r, r \le s_{N-1} \Big) \\
& = E \prod_{j\in\bigcup_{j=\overline{1,N-2}}A_j}f_j(\xi(r_j)) g_{N-1}(\xi(s_{N-1})) = \ldots = g_1(y).
\end{align*}
Since all $g_j$ are uniquelly determined, such is the distribution of $\xi.$

Define
\begin{align*}
\chi^N_\ve & =  \left(\mcP_{0,T}X_\ve(x_1, t_1, \cdot),\ldots, \mcP_{0,T}X_\ve(x_N, s_N, \cdot)\right),\\
\widehat\chi^N_\ve & =  \left(\mcP_{0,T}X^{-1}_\ve(x_1, t_1, T, T+t_1-\cdot),\ldots, \mcP_{0,T}X^{-1}_\ve(x_N, t_N, T, T+t_N-\cdot)\right), \ve\in(0;1), \\
\chi^N & =  \left(\mcP_{0,T}X(x_1, t_1, \cdot),\ldots, \mcP_{0,T}X(x_N, t_N, \cdot)\right),\\
\widehat\chi^N & =  \left(\mcP_{0,T}X^{-1}(x_1, t_1, T,T+t_1-\cdot),\ldots, \mcP_{0,T}X^{-1}(x_N, t_N, T, T+t_N-\cdot)\right). 
\end{align*}
We have proved that $\chi^N_\ve\Rightarrow \chi^N, \ve\rightarrow 0,$ in $\mcC^N([0;T])$ for any $N.$ Denote by $\pi^K$ the projection mapping in $\mcC^\infty([0;T])$ on the first $K$ coordinates, and consider the following $\mcC^\infty([0;T])-$valued elements $\chi_\ve, \widehat\chi_\ve, \chi, \widehat\chi$ defined via
\begin{align*} 
\pi^N(\chi_\ve) & = \chi^N_\ve, \pi^N(\chi)=\chi^N, \\
\pi^N(\widehat\chi_\ve) & = \widehat\chi^N_\ve, \pi^N(\widehat\chi)=\widehat\chi^N, \ve\in(0;1), N\in\mbN,
\end{align*}
their existence being guaranteed by the Kolmogorov theorem. Given $\kappa\in\mcC^\infty([s;t])$ we write $\kappa[j]$ for the $j-$th coordinate of $\kappa.$ Due to the definition of the product topology a system
$$
\big\{ \{\kappa\in\mcC^\infty([0,T])\mid \kappa[j]\in A_j, j=\overline{1,k}\}, A_j\mathrm{\ is\ open\ in\  }\mcC([0;T]), j=\overline{1,k}, k\in\mbN \big\}
$$
is a convergence-determining $\pi-$system~\cite{billingsley}[Theorem 2.2]. Consequently, $\chi_\ve\Rightarrow\chi, \ve\rightarrow 0, $ in $\mcC^\infty([0;T]).$ For the purpose of working with the inverse flows, consider a mapping $\mcI\colon \mcC^\infty([0;T])\mapsto\mcC^\infty([0;T]):$ 
\begin{align*}
\begin{cases}
\mcI(\kappa)[j](r) & = \inf\{ \kappa[i](r) \mid \kappa[i](T)\ge x_j, t_i\le r \}, r\in[t_j;T], \\
\mcI(\kappa)[j](r) & = \mcI(\kappa)[j](t_j), r\in[0;t_j), \\
j\in\mbN.
\end{cases}
\end{align*}
Since the Harris flows are right-continuous and the set $\{(x_n, t_n)\}_{n\in\mbN}$ is dense in $\mbR\times[0;T]$ $\mcI(\chi_\ve)=\widehat\chi_\ve$ a.s., $\ve\in(0;1),$ and $\mcI(\chi)=\widehat\chi$ a.s. due to~\eqref{inverse.flow:definition}. The mapping $\mcI,$ although discontinuous, is  $Law(\chi)-$continuous in the following sense. Let ${\mcC_1}$ be a set of $\psi\in\mcC^\infty([0;T])$ such that $\forall j_1, j_2\in\mbN$ 
\begin{equation*}
\big( \exists s\colon {\psi}[j_1](s)\ge {\psi}[j_2](s) \big) \Rightarrow \big({\psi}[j_1](t)\ge {\psi}[j_2](t), t\ge s\big),
\end{equation*}
and  ${\psi}[k](t_k) = x_k, k\in\mbN.$ Let $\mcC_2$ be a subset of $\mcC_1$ such that for any $\psi\in\mcC_2$
\begin{enumerate}
\item  $\forall k\in\mbN \ \exists \kappa_k>0$  
\begin{equation}
\label{prop:divergence}
\begin{cases} 
& \forall i\colon \big(x_i \ge \mcI(\psi)[k](t_i)\big) \Rightarrow \big(\psi[i]\left(T\right) - x_{k} \ge \kappa_{k}\big), \\
& \forall i\colon \big(x_i < \mcI(\psi)[k](t_i)\big) \Rightarrow \big(x_{k} - \psi[i]\left(T\right) \ge \kappa_{k}\big);
\end{cases}
\end{equation}
\item 
$\forall \delta>0 \ \forall M>0 \ \exists L\in\mbN$ 
\begin{align}
\label{prop:equicontinuity}
\sup_{l=\overline{0,\lceil T\rceil 2^{L}}}\sup_{j\colon |x_j|\le M, t_j = l2^{-L}} \sup_{\Delta t\in[0;2^{-L}]} \big| \psi[j](t_j+\Delta t)- x_j \big| \le \delta.
\end{align}
\end{enumerate}
We state that if $\psi_n\rightarrow\psi,n\rightarrow\infty, \psi_n\in\mcC_1, n\in\mbN, \psi\in{\mcC}_2,$ then $\mcI(\psi_n)\rightarrow\mcI(\psi),n\rightarrow\infty.$ To see that, suppose the opposite. Then there exist a sequence $\{ s_n\}_{n\in\mbN}\subset[0;T]^\infty$ and numbers $s\in[0;T], \kappa\in\mbR^+$ and $k\in\mbN$ such that $s_n\rightarrow s,n\rightarrow\infty,$ and
$$
\inf_{n\in\mbN} |\mcI(\psi_n)[k](s_n) - \mcI(\psi)[k](s_n)| \ge \kappa.
$$
 
Since $\mcI(\psi)[k]$ is a continuous function,
$$
\liminf_{n\rightarrow\infty} \left| \mcI(\psi_n)[k](s_n) - \mcI(\psi)[k](s) \right| \ge \frac{\kappa}{2},
$$
so there exists a sequence $\{j_n\}_{n\in\mbN}, t_{j_n} \le s_n, n\in\mbN,$ such that at least one of the following relations holds:
\begin{align}
\label{eq:prelimit.1}
 \psi_n[j_n](s_n) - \mcI(\psi)[k](s) \ge \frac{\kappa}{4} \mathrm{\ and\ } \left( \psi_n[j_n](T) < x_k \mathrm{\ infinitely\ often}\right), \\
 \label{eq:prelimit.2}
\mcI(\psi)[k](s) - \psi_n[j_n](s_n) \ge \frac{\kappa}{4} \mathrm{\ and\ } \left( \psi_n[j_n](T) \ge  x_k \mathrm{\ infinitely\ often}\right).
\end{align}
We suppose that \eqref{eq:prelimit.1} holds. The case \eqref{eq:prelimit.2} is treated similarly.

Because of \eqref{prop:equicontinuity}, there exist $\ve>0$ and $j\in\mbN$ such that 
\begin{align}
\label{eq:close.point}
\begin{cases}
t_{j} \le s-\ve, \\
\sup\limits_{t\in s - \ve;s+\ve]} \left|\psi[j](t) - \mcI(\psi)[k](s)\right| \le \frac{\kappa}{8}, \\
\psi[j](s) - \mcI(\psi)[k](s) > 0.
\end{cases}
\end{align}
Due to \eqref{eq:close.point} and \eqref{prop:divergence} we have: 
\begin{equation}
\label{eq:close.point.T}
\psi[j](T) \ge x_k +\kappa_k.
\end{equation}
Combining \eqref{eq:close.point} with \eqref{eq:prelimit.1}, we have:
$$
\inf_{n\in\mbN} \big(\psi_n[j_n](s_n) - \psi[j](s_n)\big) \ge \frac{\kappa}{8}.
$$
At the same time, $\psi_n[j]\rightarrow\psi[j], n\rightarrow\infty,$ in $\mcC([0;T])$ thus there exists $n_0$ such that
\begin{equation}
\label{eq:limit.proc.diver}
\forall n\ge n_0 \ \big(\psi_n[j_n](s_n) - \psi_n[j](s_n)\big) \ge \frac{\kappa}{16}.
\end{equation}

Moreover, \eqref{eq:limit.proc.diver} and \eqref{eq:prelimit.1} imply that
$$
\psi[j](T) = \lim\limits_{n\rightarrow\infty} \psi_n[j](T) \le \liminf\limits_{n\rightarrow\infty} \psi_n[j_n](T) < x_k,
$$
which contradicts \eqref{eq:close.point.T}. Hence $\mcI(\psi_n)\rightarrow\mcI(\psi),n\in\infty,$ in $\mcC^{\infty}([0;T]).$

Obviously, $\chi_\ve\in{\mcC}_1$ a.s., $\ve\in(0;1).$ In order to check that $\chi\in\mcC_2,$ only the properties~\eqref{prop:divergence} and~\eqref{prop:equicontinuity} need be verified. For arbitrary positive $\delta$ and $M$
\begin{align}
\label{eq:equicontinuity}
P\left\{ \bigcup_{l=\overline{0,\lceil T\rceil 2^{L}}} \left\{\sup_{u\in[-M;M]}\sup_{t\in[l2^{-L};(l+1)2^{-L}]} \left| X(u,l2^{-L},t) - u\right| \ge \delta\right\}  \mathrm{\ for\ inf.\ many\ }L \right\} = 0.
\end{align}
Indeed, if
$$
U_{lL} = \left\{ \sup_{u\in[-M;M]}\sup_{t\in[l2^{-L};(l+1)2^{-L}]} \left| X(u,l2^{-L},t) - u\right| \ge \delta \right\},
$$
then it is sufficient for \eqref{eq:equicontinuity} to hold that the series $\sum_{L\ge 1}\sum_{l=\overline{0,\lceil T\rceil 2^{L}}} P\left(U_{lL}\right)$ converge. Proceeding exactly as in the proof of Theorem 4.7 of~\cite{harris} and using an estimate~\cite{harris}[4.8] one obtains:
\begin{align*}
\sum_{L\ge 1}\sum_{l=\overline{0,\lceil T\rceil 2^{L}}} & P\left(U_{lL}\right) \le \left( \frac{16M}{\delta} + 2\right) \sqrt\frac{2}{\pi}\sum_{L\ge 1}  \left(\lceil T\rceil 2^L +1\right) \int\limits^\infty_{\delta 2^{\frac{L}{2}}} e^{-\frac{y^2}{2}}dy \\
& \le \left( \frac{16M}{\delta} + 2\right)  \frac{2}{\delta} \sqrt\frac{2}{\pi} \lceil T\rceil  \sum_{L\ge 1} 2^{\frac{L}{2}} e^{-\delta 2^{L}} < +\infty.
\end{align*}
It is easily seen that \eqref{eq:equicontinuity} implies the property~\eqref{prop:equicontinuity}. To prove \eqref{prop:divergence}, note that, according to~\cite{harris}[Chapter 7],  a mapping $X(\cdot,s,t)$ is a jump function, so, due to~~\cite{harris}[Chapter 4], for any $s\in[0;T)$:
\begin{align*}
P & \left\{ \{X(x,s,T)\mid x \in\mbR\} \cap \{x_1,\ldots, x_N\} \not= \emptyset \right\} \\
 & = P\left\{ \left\{X(x,s,T)\mid x = \frac{u}{v}, u,v\in\mathbb{Z}\right\} \cap \{x_1,\ldots, x_N\} \not= \emptyset  \right\} \\
 & \le \sum_{u,v\in\mathbb{Z}} P\left\{ X\left(\frac{u}{v},s,T\right)\in \{x_1,\ldots, x_N\}  \right\} =0,
\end{align*}
and thus follows \eqref{prop:divergence}. Hence $\chi\in\mcC_2$ a.s..

The continuous mapping theorem~\cite{kallenberg}[Theorem 4.27 + Exercise 4.27] implies that
$$
\left(\chi_\ve,\widehat\chi_\ve\right) \Rightarrow \left( \chi,\widehat\chi \right), \ve\rightarrow 0+,
$$
which is essentially, after a reformulation, the first assertion of the theorem. 

To verify the second assertion of the theorem we start with checking that given $s, t, s <t,$ it holds that  $\mu_\ve(s,t)\Rightarrow \mu(s,t)$ in $\mcM(\mbR)$ as $\ve\to0+,$ or, equivalently~\cite{kallenberg:random.measures}[Theorem 4.2], that for any continuous compactly supported function $f$ $\lg \mu_\ve(s,t), f\rg\Rightarrow\lg \mu(s,t), f\rg, \ve\rightarrow0,$ where for any $\nu\in\mcM(\mbR)$ and any function $g$ 
$$
\lg \nu, g \rg \coloneqq \int_\mbR g(y) \nu(dy),
$$ 
assuming that the integral exists. We proceed gradually. Firstly, we prove this convergence to hold with an additional assumption of $f$ being Lipshitz continuous; secondly, we establish the result for arbitrary continuous compactly supported functions. 

Suppose that $f$ is compactly supported and Lipschitz continuous with a Lipschitz constant $C_f.$ Let $\supp(f)\subset[-S;S]$ for some $S>0$ and $R_f=\sup_{y\in\mbR}|f(y)|.$ Given a standard Brownian motion $W$ it holds, for any $\ve\in(0;1)$ and $M>S,$ that
\begin{align*}
 E& \Big|\int\limits_{|y|\ge M} f\left(X_\ve\left(y,s,t\right)\right)dy \Big| \le 
E \sum_{k\ge M}\int\limits_k^{k+1} \big(|f\left(X_\ve\left(y,s,t\right)\right)|  + |f\left(X_\ve\left(-y,s,t\right)\right)|\big)dy \\
& \le  R_f \sum_{k\ge M}\int\limits_k^{k+1} \Big({ P\left\{ X_\ve\left(y,s,t\right) \in[-S;S]\right\}  + P\left\{ X_\ve\left(-y,s,t\right)\in[-S;S]\right\} }\Big)dy \\
& \le 2 R_f \sum_{k\ge M} P\left\{ X_\ve\left(k,s,t\right) \le S\right\} 
\le 2 R_f \sum_{k\ge M} P\left\{ W(t-s) \ge k-S\right\}.
\end{align*}
The same estimate, obviously, holds for $X,$ too. Fix $\delta >0.$ Then there exists $M$ such that 
\begin{align}
\label{eq:int.outside.estimate}
\max\Big\{ \sup_{\ve(0;1)}E\Big|\int_{|y|\ge M} f\left(X_\ve\left(y,s,t\right)\right)dy\Big|, E\Big|\int_{|y|\ge M} f\left(X\left(y,s,t\right)\right)dy\Big| \Big\} \le \delta.
\end{align}
Define
\begin{align*}
\lg \mu^N_\ve, f\rg & = \sum_{k=-N}^{N} f\left(X_\ve\left(\frac{kM}{N}, s, t\right)\right) MN^{-1}, \\
\lg \mu^N, f\rg & = \sum_{k=-N}^{N} f\left(X\left(\frac{kM}{N}, s, t\right)\right) MN^{-1}, N\in\mbN.
\end{align*}
Fix a Lipschitz continuous function $g$ with the Lipschitz constant $C_g.$ Put 
\begin{align*}
A_{N} & = Eg\left(\lg \mu^N(s,t), f \rg\right) - Eg\left(\lg \mu(s,t), f \rg\right), \\
A_{\ve N} & = Eg\left(\lg \mu_\ve^N(s,t), f \rg\right) - Eg\left(\lg \mu_\ve(s,t), f \rg\right), \ve\in(0;1). 
\end{align*}
Since Harris flows are stationary w.r.t. the time variable, we have, by~\eqref{eq:int.outside.estimate}, that
\begin{align}
\label{a_e:estimate.deriving}
|A_{\ve N}| & = \left|Eg\left(\lg \mu_\ve(s,t), f \rg\right) - Eg\left(\lg \mu^N_\ve(s,t), f \rg\right)\right|  \nonumber \\
&  \le C_g E \left|\lg \mu_\ve(s,t), f \rg - \lg \mu^N_\ve(s,t), f \rg\right| \nonumber \\
& \le C_g\delta  + C_g E \left| \int_{-M}^M f\left(X_\ve\left(y,s,t\right)\right)dy - \sum_{k=-N}^{N} f\left(X_\ve\left(\frac{kM}{N}, s, t\right)\right) MN^{-1} \right| \nonumber \\
& \le C_g\delta  +  C_g E \sum_{k=-N}^{N} \int_\frac{kM}{N}^\frac{(k+1)M}{N} \left| f\left(X_\ve\left(y,s,t\right)\right) - f\left(X_\ve\left(\frac{kM}{N}, s, t\right)\right) \right| dy \nonumber \\
& = C_g\delta  +  2C_g N E\int_0^{MN^{-1}} \left| f\left(X_\ve\left(y,s,t\right)\right) - f\left(X_\ve\left(0, s, t\right)\right) \right| dy \nonumber \\
& \le C_g\delta  +  2 C_g C_f N \int_0^{MN^{-1}}  \left[ E \left( X_\ve\left(y,s,t\right) - X_\ve\left(0, s, t\right)\right)^2 \right]^{\frac{1}{2}}  dy.
\end{align}
The result of~\cite{dorogovtsev:entropy}[Lemma 5], after an investigation of its proof, can be reformulated as follows: for any Harris flow $Y$ and any $y_1, y_2 \in\mbR, 0\le s \le t, t-s\le 1,$
\begin{equation}
\label{square.estimate}
E\left(Y(y_1,s,t) - Y(y_2,s,t)\right)^2 \le (y_1-y_2)^2 + \frac{8}{\pi}|y_1-y_2|.
\end{equation}
Using~ \eqref{a_e:estimate.deriving} and~\eqref{square.estimate} one obtains, for $MN^{-1}\le 1:$
\begin{equation*}
\label{a_e:estimate}
|A_{\ve N}| \le C_g\delta  +  6C_f C_g N\int_{0}^{MN^{-1}} y^{\frac{1}{2}}dy \le C_g\delta  + 4 C_g C_f M^{\frac{3}{2}}N^{-\frac{1}{2}}.
\end{equation*}
Exactly the same reasoning is applicable in the case of $A_N.$ 
So, for sufficiently large $N,$
\begin{align*}
\sup_{\ve\in(0;1)}\left|A_{\ve N}\right|  +\left|A_N\right|\le 2C_g\delta  + 8 C_g C_f M^{\frac{3}{2}}N^{-\frac{1}{2}}.
\end{align*}
Therefore,
\begin{align}
\label{ineq:lipsh.estimate}
\Big| E g\left(\lg \mu_\ve(s,t), f \rg\right)- Eg\left(\lg \mu(s,t), f \rg\right) \Big| \le &  \left|Eg\left(\lg \mu^N_\ve(s,t), f \rg\right) - Eg\left(\lg \mu^N(s,t), f \rg\right) \right| \nonumber \\ 
& + 2C_g\delta  + 8 C_g C_f M^{\frac{3}{2}}N^{-\frac{1}{2}},
\end{align} 
for sufficiently large $N.$ Here $\delta$ can made arbitrary small by taking $M$ large enough. 
Due to the first statement of the theorem, for any fixed natural $N,M$
\begin{align*}
\Bigg(X_\ve(-M,s,t), \ldots,& X_\ve\left(-\frac{kM}{N},s,t\right),  \ldots, X_\ve(0,s,t), \ldots,  \\ 
& X_\ve\left(\frac{kM}{N},s,t\right), \ldots, X_\ve(M,s,t)\Bigg) \nonumber  \Rightarrow \\
   \Bigg(X(-M,s,t), \ldots, & X\left(-\frac{kM}{N},s,t\right), \ldots, X(0,s,t), \ldots, \\
&  X\left(\frac{kM}{N},s,t\right), \ldots, X(M,s,t)\Bigg) 
\end{align*}
in $\mbR^{2N+1},$ as $\ve\rightarrow 0.$ Hence for any $N\in\mbN$
\begin{equation*}
\label{second.argument}
 Eg\left(\lg \mu^N_\ve(s,t), f \rg\right) - Eg\left(\lg \mu^N(s,t), f \rg\right) \to 0,\ve\to 0.
\end{equation*}
This, together with \eqref{ineq:lipsh.estimate}, implies that for any Lipschitz continuous $g$ 
$$
Eg\left(\lg \mu_\ve(s,t), f\rg\right)\to Eg\left(\lg \mu(s,t), f\rg\right), \ve\to 0.
$$
Equivalently, 
\begin{equation}
\label{weak.cong.form}
\lg \mu_\ve(s,t), f \rg \Rightarrow\lg \mu(s,t), f \rg\mathrm{\ in\ }\mbR, \ve\to\ 0.
\end{equation} 

Now we shall show that \eqref{weak.cong.form} holds for arbitrary continuous functions $f$ whose support is contained in $[-S;S].$  It is sufficient to show that for any $\delta > 0$ there exists a Lipschitz continuous function $f^\star$ such that for any Lipschitz continuous function $g$
\begin{align*}
& \left| Eg\left(\lg \mu_\ve(s,t), f\rg\right)-Eg\left(\lg \mu_\ve(s,t), f^\star\rg\right)\right| \le C\delta, \\
& \left| Eg\left(\lg \mu(s,t), f\rg\right)-Eg\left(\lg \mu(s,t), f^\star\rg\right)\right| \le C \delta,
\end{align*}
where the constant $C$ does not depend on $f^\star$ or $\ve.$ Obviously, there exists a Lipschitz continuous function $f^\star$ supported on $[-S;S]$ such that $\max_{y\in\mbR}|f(y)-f^\star(y)|\le \delta.$ Let $\nu$ denote any of $\{\mu_\ve(s,t)\}_{\ve\in(0;1)},$ or $\mu(s,t).$ We have:
\begin{align*}
\left| Eg\left(\lg \nu,f\rg\right) -Eg\left(\lg \nu,f^\star\rg\right)\right| \le C_gE\int_{-S}^S |f(y)-f^\star(y)| \nu(dy) \le C_g \delta E \nu((-S;S]).
\end{align*}
Let $W$ be a standard Wiener process started from $0.$ Then
\begin{align*}
\label{lip.estimate}
E \mu(s,t)((-S;S]) & = E\lambda\{y\mid X(y, s, t)\in(-S;S]\} \nonumber \\
& = \int_0^\infty P\big\{ \lambda\{y\mid X(y, s, t)\in(-S;S]\} \ge c  \big\} dc \nonumber \\
& \le \int_0^\infty P\big\{ \exists y_1, y_2\colon y_2 - y_1 \ge c, X(y_1,s,t), X(y_2,s,t) \in[-S;S]  \big\} dc \nonumber \\
& \le 2S + \int_{2S}^\infty P\Big\{ \exists y\colon |y|\ge \frac{c}{2},  X(y,s,t)\in[-S;S]  \Big\} dc \nonumber \\
& \le 2S + \int_{2S}^\infty P\Bigg\{ X\left(\frac{c}{2},s,t\right)\le S\mathrm{\ or \ }X\left(\frac{-c}{2},s,t\right)\ge -S  \Bigg\} dc \nonumber \\
& \le 2S + 2 \int_{2S}^\infty P\left\{ \frac{c}{2}+W(t-s)\le S\right\} dc  \nonumber \\
& = 2S + 2 \int_{2S}^\infty P\left\{ W(t-s) \ge \frac{c}{2}-S\right\} dc \nonumber \\
& = 2S + 4 \int_0^\infty P\big\{ W(t-s) \ge c\big \} dc \le 2S + 2E|W(t-s)| \le 2S + 2T.
\end{align*} 
This estimate holds also for $\mu_\ve(s,t).$ Thus
$$
\left| Eg\left(\lg \nu,f\rg\right) -Eg\left(\lg \nu,f^\star\rg\right)\right|  \le C_g (2S+2T)\delta,
$$
and, as a result, \eqref{weak.cong.form} holds for compactly supported continuous $f.$ So $\mu_\ve(s,t)\Rightarrow\mu(s,t)$ in $\mcM(\mbR), \ve\to 0.$

By \cite{harris}[Theorem 10.5] we have that $Law(\mu_\ve(s,t))=Law(\hat\mu_\ve(s,t)), \ve\in(0;1),$ and $Law(\mu(s,t))=Law(\hat\mu(s,t))$ so, since $\mcM(\mbR)$ is separable the mappings $\mu_\ve(s_1,t_1)$ and $\mu_\ve(s_2,t_2)$ are independent as soon as $(s_1,t_1)\cap(s_2,t_2)=\emptyset,$ as well as those of the inverse flows. Hence, the second assertion of the theorem follows by standard reasoning.
\qed
\begin{remark}
In Theorem~\ref{harris:approximation} the measures $\mu_\ve(s,t)$ are considered instead of the mappings $X_\ve(\cdot,s,t)$ (although the finite-dimensional distributions of the latter ones are convergent as random variables) because the family $\{X_\ve(\cdot,s,t)\}_{\ve\in(0;1)}$ is not tight in $\mcD(\mbR):$ the limit function $X(\cdot, s,t)$ is discontinuous while the prelimit ones are continuous, and the Skorokhod topology does not allow such convergence to happen. However, one can consider another possible topologies, for instance the weak convergence topology metriced with the Levy-Prokhorov distance (or, equivalently, with the $L_1-$Wasserstein metric with a bounded integrand). This is essentially what is done in the theorem.
\end{remark}
\begin{remark}
In~\cite{dorogovtsev:rate.harris.flows} estimates on the Wasserstein metric between the distributions of the forward n-point motions of one-dimensional Harris flows with compactly supported covariance functions are obtained in terms of the diameters of supporting sets. 
However, in our case for any collection $\{\vf_\ve\}_{\ve\in(0;1)}$ of compactly supported functions approximating $\vf,$ the corresponding supporting sets grow indefinitely. 
\end{remark}

\begin{section}{Acknowledgments}
The author is indebted to a referee for helpful comments and suggestions. 
\end{section}

\end{document}